\documentclass[12pt,a4paper]{article} 
\usepackage{a4}

\usepackage[intlimits,tbtags]{amsmath}
\usepackage{amsfonts}
\usepackage{amsthm,amscd}

\usepackage{cite}
\usepackage{latexsym}

\newcommand{\E}{\mathrm{e}}

\newcommand{\tr}{\triangleright}
 
\newcommand{\R}{\mathcal{R}}
\newcommand{\F}{\mathcal{F}}

\newcommand{\id}{\mathrm{id}}

\newcommand{\ad}{{\mathrm{ad}\,}}
\newcommand{\jhat}{{\hat{\jmath}}}

\newcommand{\Xcal}{\mathcal{X}}

\newcommand{\Ocal}{\mathcal{O}}

\newcommand{\Ket}[1]{\lvert #1\rangle}

\newcommand{\Braket}[3]{\langle #1\rvert #2 \lvert #3\rangle}
\newcommand{\rBraket}[3]{\langle #1\rVert #2 \lVert #3\rangle}

\newtheorem{Theorem}{Theorem}
\newtheorem{Proposition}{Proposition}

\newcommand{\h}{\hbar}

\newcommand{\Usu}{{\mathcal{U}(\mathrm{su}_2)}}

\newcommand{\Uqsu}{{\mathcal{U}_q(\mathrm{su}_2)}}
\newcommand{\Uhsu}{{\mathcal{U}_\h(\mathrm{su}_2)}}
\newcommand{\Uhg}{{\mathcal{U}_\h(\mathfrak{g})}}
\newcommand{\Ug}{{\mathcal{U}(\mathfrak{g})}}

\newcommand{\CG}[6]{
  \biggl(\begin{matrix}#1\!\! &  #2 \\#4\!\! & #5 \end{matrix}
  \biggm|
  \begin{matrix}
    \vphantom{#1#2} {#3} \\ \vphantom{#4#5} {#6}
  \end{matrix}\biggr)}

\newcommand{\CGq}[6]{
  \biggl(\begin{matrix}#1\!\! &  #2 \\#4\!\! & #5 \end{matrix}
  \biggm|
  \begin{matrix}
    \vphantom{#1#2} {#3} \\ \vphantom{#4#5} {#6}
  \end{matrix}\biggr)_{\!\!q}}

\newcommand{\RCs}[6]{{\mathrm{R}^{#1#2#3}_{#4#5#6}}}

\newcommand{\SixJR}[6]{
  \left\{\begin{matrix}#1&#2&#4\\#3&#6&#5\end{matrix}\right\}}

\newcommand{\SixJ}[6]{
  \left\{\begin{matrix}#1&#2&#3\\#4&#5&#6\end{matrix}\right\}}

\newcommand{\NineJ}[9]{
  \left\{\begin{matrix}#1&#2&#3\\#4&#5&#6 \\ #7&#8&#9\end{matrix}\right\}}
    
{\begin{list}{}{ \settowidth{\labelwidth}{\textrm{#1}}
      \setlength{\leftmargin}{\parindent}
      \addtolength{\leftmargin}{\labelwidth}
      \addtolength{\leftmargin}{\labelsep}}
      }%
    {\end{list}}

\begin{document}

\rightline{Preprint: IUB-TH-0412}

\vspace{4em}
\begin{center}
  
  {\Large{\bf Reconstruction of universal Drinfeld twists from representations}}

\vspace{3em}

\textbf{Christian Blohmann}

\vspace{1em}
 
International University Bremen, School of Engineering and Science\\
Campus Ring 1, 28759 Bremen, Germany\\[1em]

\end{center}

\vspace{1em}

\begin{abstract}
  Universal Drinfeld twists are inner automorphisms which relate the coproduct of a quantum enveloping algebra to the coproduct of the undeformed enveloping algebra. Even though they govern the deformation theory of classical symmetries and have appeared in numerous applications, no twist for a semi-simple quantum enveloping algebra has ever been computed. It is argued that universal twists can be reconstructed from their well known representations. A method to reconstruct an arbitrary element of the enveloping algebra from its irreducible representations is developed. For the twist this yields an algebra valued generating function to all orders in the deformation parameter, expressed by a combination of basic and ordinary hypergeometric functions. An explicit expression for the universal twist of su(2) is given up to third order.
\end{abstract}


\newpage
\section{Introduction}

Quantum enveloping algebras are formal deformations of the enveloping Hopf algebras of Lie algebras \cite{Drinfeld:1989}. While the notion of quantum enveloping algebras is very general, comprising arbitrary deformations, the most famous examples are the Drinfeld-Jimbo deformations \cite{Drinfeld:1985,Jimbo:1985} which act as natural symmetry structures on quantum spaces \cite{Manin:1988,Faddeev:1990,Carow-Watamura:1990}. Drinfeld observed that every quantum enveloping algebra is related to the corresponding undeformed enveloping algebra by an inner automorphism which he called universal twist \cite{Drinfeld:1989b} and which now bears his name. Given the universal Drinfeld twist one can reconstruct the corresponding quantum enveloping algebra up to isomorphism. In this sense, the twist contains all information on the quantum deformation of a classical symmetry \cite{Gerstenhaber:1992b}. 

Due to their pivotal role for the deformation theory of symmetries, universal Drinfeld twists have found numerous important applications beyond the structure theory of quantum enveloping algebras, such as to quantum statistics on quantum spaces \cite{Fiore:1996}, quantum spin chains \cite{Terras:1999,Maillet:2000}, noncommutative quantum field theory \cite{Grosse:2001}, or to algebraic geometry \cite{Racinet:2002}, just to name a few recent examples. Our original motivation was the observation that certain twists yield a covariant realization of quantum spaces by a star product \cite{Giaquinto:1992,Blohmann:2002a} within the framework of deformation quantization \cite{Bayen:1978a}. Such a description of noncommutative spaces by formal deformations of algebras \cite{Gerstenhaber:1964} has appeared naturally in the context of string theory \cite{Seiberg:1999}, the construction of gauge theories on noncommutative spaces \cite{Madore:2000b,Jurco:2001}, and the subsequent development of noncommutatvie quantum field theories. (For reviews of noncommutative field theories see 
\cite{Douglas:2001} and \cite{Szabo:2001}.)

The noncommutative geometry on which so far most noncommutative quantum field theories have been constructed is defined by constant commutators of the space-time observables. Such a noncommutativity breaks Lorentz symmetry, which had to be expected because the constant commutator can be viewed as due to a constant background field, in string theory a constant $B$-field on a
D-brane. It was hoped that a small noncommutativity would lead to an equally small violation of Lorentz symmetry. However, on the level of regularization of loop diagrams an interdependence of ultra-violet and infra-red cutoff scales appears \cite{Minwalla:1999,Matusis:2000} which seems to put even large scale Lorentz symmetry and weakened notions of locality of noncommutative quantum field theory into doubt \cite{Alvarez-Gaume:2003}, indicating that the breaking of symmetries is not under good control --- at least for the case of constant noncommutativity. These serious deficiencies seem to be reason enough to reconsider such deformations, for which the symmetry structure can be deformed together with the space, so that covariance is preserved. That is, quantum spaces \cite{Manin:1988,Faddeev:1990,Carow-Watamura:1990} which carry a covariant representation of the Drinfeld-Jimbo deformation \cite{Drinfeld:1985,Jimbo:1985} of the enveloping symmetry algebra. Three particularly important quantum spaces, the quantum plane, quantum Euclidean 4-space, and quantum Minkowski space, have been shown to be realizable as deformation quantization \cite{Giaquinto:1992} by universal Drinfeld twists \cite{Blohmann:2002a}.

While the quantum enveloping algebra can be reconstructed rather easily given the universal Drinfeld twist, there is no general soution for the inverse problem of calculating the Drinfeld twist for a given quantum deformation. The existence of twists is proved by homological methods which are inherently nonconstructive. To our best knowledge, no universal Drinfeld twist for the Drinfeld-Jimbo deformation of a semisimple Lie algebra has ever been computed successfully, not even for the simplest possible case of $\mathrm{su}_2$.  (In \cite{Bonechi:1995} and \cite{Fiore:1997} the non-semisimple case of the Heisenberg algebra was studied.) It could be argued that the universal twist is more or less the square root of the universal $\R$-matrix, so that Drinfeld's ingenious but simple construction of the $\R$-matrix by the dual pairing of the Borel Hopf subalgebras might be used. For the case of triangular deformations this reasoning appears to lead, indeed, to a method to construct the twist \cite{Stolin:2003}. For the non-triangular Drinfeld-Jimbo deformations, however, this argument falls short, as is confirmed by the complexity of the expressions derived in \cite{Dabrowski:1996}. But even though there are no closed form expressions or simple constructions for the twist, one might expect that a brute force calculation by means of a computer algebra system ought to be possible up to an order of the perturbation parameter high enough for all reasonable applications. However, it turns out that the naive attempt of an algebraic order by order calculation quickly runs into overwhelming combinatorial problems, as it was described in \cite{Dabrowski:1996} where the authors did not go beyond the second order. 

A closer inspection of the brute force approach reveals, firstly, that the number of operations which have to be carried out increases at least exponentially with the order of the perturbation parameter, so the algorithm is certainly nonpolynomial. Secondly, the results expressed in terms of the Poincar\'e-Birkhoff-Witt basis of the enveloping algebra are extremely lengthy and do not appear to povide any structural insight. Thirdly, it is unclear how to implement the algorithm such that it yields the particular twist which realizes the star product of quantum spaces. In conclusion, it is fair to say that the calculation of universal Drinfeld twists turns out to be a computational problem in any respect.

While little is known about the calculation of universal twists in the enveloping algebra, their representations are well understood and have been computed explicitly. They are essentially given by a contraction of deformed and undeformed Clebsch-Gordan coefficients as it was first observed in \cite{Curtright:1991}. For the Drinfeld-Jimbo $q$-deformation the Clebsch-Gordan coefficients are known explicitly, so we obtain the representations of the twist in a closed form, expressed by basic and ordinary hypergeometric functions. The approach to the calculation of universal twists presented here is to reconstruct the twists from their representations. The obvious advantage of this approach is that, rather than starting with algebraic calculations from scratch, it builds on the computational effort which has gone into the calculation of $q$-Clebsch-Gordan coefficients. The explicit calculations in this article are carried out for $\Uqsu$ although the methods will be seen to be generic.

We will proceed as follows: In order to make this article reasonably self-contained we will start in Sec.~\ref{sec:Drinfeld} with a short introduction to universal Drinfeld twists, giving some basic results which will be refered to in the following. Moreover, we briefly explain why universal twists appear naturally in the context of star products. In Sec.~\ref{sec:Method} we will develop a method to reconstruct an element in the enveloping algebra from its irreducible representations. From Lie theory it is clear that in the semi-simple case the representations determine the algebra element uniquely. But how do we actually compute it? The key to developing a constructive method is the choice of a suitable basis of $\Uqsu$. We will choose the basis of tensor operators because the matrix elements satisfy useful orthogonaltiy relations, which then lead to the desired reconstruction method. The results are given in Proposition~\ref{th:Reconst1} and Eq.~\eqref{eq:Reconst7}. In Sec.~\ref{sec:Reconstruction} we apply the reconstruction method to the representations of the universal Drinfeld twist. This leads to the main result presented here: a formula for the the universal Drinfeld twist, given in Eq.~\eqref{eq:Reconst11} by an algebra-valued generating function in $q=\E^\h$. In order to obtain the twist to each order in the perturbation parameter $\h$, we yet need to expand the generating function in powers of $\h$. Surprisingly, the problem of perturbative expansion of basic hypergeometric functions has recieved little attention in the literature. In particular, no closed formulas for such expansions have been derived as yet.  While a thorough study of this problem is beyond the scope of this article, we present the first few steps in this direction which suffice to make the order of order expansion of the generating function of the twist accessible to efficient computer algebra calculations. This is exemplified in Eq.~\eqref{eq:Fresult} where the universal twist of $\Uhsu$ was computed up to third order. In Sec.~\ref{sec:Conclusion} we concludingly assess the computational value of the generating function of the twist and indicate how the approach presented here will be naturally continued.

\section{Brief introduction to Drinfeld twists}
\label{sec:Drinfeld}

For the reader's convenience we briefly review how Drinfeld twists appear naturally in the study of of formal deformations of algebras and Hopf algebras. The approach and the results described here are essentially due to Gerstenhaber \cite{Gerstenhaber:1964} and Drinfeld \cite{Drinfeld:1989,Drinfeld:1989b}. 
The formal perturbation parameter is $\h$, the completion of a complex vector
space or algebra $A$ with respect to the $\h$-adic topology by formal power series is denoted as usual by $A[[\h]]$.

An $\h$-adic algebra $A'$ is called a deformation of an algebra $A$ if
$A'/\h A'$ and $A$ are isomorphic as algebras. Analogously, an
$\h$-adic Hopf algebra $H'$ is called a deformation of a Hopf algebra
$H$ if $H'/ \h H'$ and $H$ are isomorphic as Hopf algebras. Recall,
that $\Ug$ is a Hopf algebra with the canonical Lie Hopf structure
defined on the generators $g\in \mathfrak{g}$ by the coproduct $\Delta(g)
= g \otimes 1 + 1 \otimes g$, counit $\varepsilon(g) = 0$, and
antipode $S(g) = -g$. The Drinfeld-Jimbo algebra $\Uhg$ is a
deformation of this Hopf algebra $\Ug$. This can be seen by developing
the commutation relations and the Hopf structure of $\Uhg$ as formal
power series in $\h$ and keeping only the zeroth order terms, which
yields the commutation relations and the Lie Hopf structure of $\Ug$.

Gerstenhaber has shown \cite{Gerstenhaber:1964} that whenever the
second Hochschild cohomology of $A$ with coefficients in $A$ is zero,
$H^2(A,A) = 0$, then all deformations of $A$ are trivial up to
isomorphism. That is, any deformation $A'$ of $A$ is isomorphic to the
$\h$-adic completion of the undeformed algebra, $A' \cong A[[\h]]$.
Algebras with this property are called rigid. The second Whitehead
lemma states that the second Lie algebra cohomology of a semisimple
Lie algebra $\mathfrak{g}$ and, hence, the second Hochschild
cohomology of its enveloping algebra is zero. Therefore, the
enveloping algebra $\Ug$ of a semisimple Lie algebra $\mathfrak{g}$ is
rigid. In particular, there is an isomorphism of algebras $\alpha:
\Uhg \rightarrow \Ug[[\h]]$, by which the the Hopf structure
$\Delta'$, $\varepsilon'$, $S'$ of $\Uhg$ can be transfered to
$\Ug[[\h]]$,
\begin{equation}
\label{eq:HopfDeform}
  \Delta_\h := (\alpha\otimes\alpha)
    \circ \Delta' \circ\alpha^{-1} \,,\quad
  \varepsilon_\h := \varepsilon' \circ \alpha^{-1} \,,\quad
  S_\h := \alpha \circ S' \circ \alpha^{-1} \,,
\end{equation}
such that $\alpha$ becomes an isomorphism of Hopf algebras from $\Uhg$
to $\Ug[[\h]]$ with this deformed Hopf structure. Let $\alpha'$ be
another such isomorphism and $\Delta'_\h$, $\varepsilon'_\h$, $S'_\h$
be defined as in Eq.~\eqref{eq:HopfDeform} with $\alpha'$ instead of
$\alpha$. Then $\alpha'$ is an isomorphism of Hopf algebras from
$\Uhg$ to $\Ug[[\h]]$ with the primed Hopf structure,
\begin{equation}
  (\Ug[[\h]], \Delta_\h, \varepsilon_\h, S_\h)
  \stackrel{\alpha}{\longleftarrow} \Uhg
  \stackrel{\alpha'}{\longrightarrow}
  (\Ug[[\h]], \Delta'_\h, \varepsilon'_\h, S'_\h) \,,
\end{equation}
hence, $\alpha'\circ \alpha^{-1}$ is an isomorphism of Hopf algebras.
We conclude that, while the Hopf structure~\eqref{eq:HopfDeform}
may depend on the isomorphism $\alpha$, it is unique up to an
isomorphism of Hopf algebras. 

As a consequence of the first Whitehead lemma, the first Hochschild
cohomology of the enveloping algebra $\Ug$ of a semisimple Lie algebra
is zero. This implies, that the two homomorphisms $\Delta$ and
$\Delta_\h$ from $\Ug[[\h]]$ to $(\Ug \otimes \Ug)[[\h]]$ with
$\Delta_\h = \Delta + \Ocal(\h)$ are related by an inner automorphism,
as it was observed by Drinfeld \cite{Drinfeld:1989,Drinfeld:1989b}.

\begin{Theorem}[Drinfeld]
\label{th:Drinfeld2}
Let $\mathfrak{g}$ be a semisimple Lie algebra, and let $\Delta_\h$ be
defined as in Eq.~\eqref{eq:HopfDeform}. Then there is an invertible
element $\F \in \bigl(\Ug \otimes \Ug\bigr)[[\h]]$ such that
$\Delta_\h(g) = \F^{-1} \Delta(g) \F$, which is called a Drinfeld
twist from $\Delta$ to $\Delta_\h$.
\end{Theorem}

It can be shown that such a Drinfeld twist not only relates the deformed and undeformed coproducts but also the counits and antipodes. Hence, a universal Drinfeld twist uniquely determines the corresponding quantum enveloping algebra. In that sense the twist contains the entire structural information on a quantum deformation of an enveloping algebra. The twist of Theorem~\ref{th:Drinfeld2} is not unique. For a given quantum enveloping algebra any two twists are related by a noncommutative 2-coboundary in the sense of \cite{Majid}.

Drinfeld has shown, that the isomorphism of $\Ug[[\h]]$ and $\Uhg$ can be
chosen to leave a given Cartan subalgebra invariant:
\begin{Theorem}[Drinfeld \cite{Drinfeld:1989}, Prop.~4.3]
  \label{th:CartanPreserve}
  Let $\mathfrak{g}$ be a semisimple Lie algebra and $\mathfrak{h}
  \subset \mathfrak{g}$ a Cartan subalgebra. Then there exists an
  isomorphism of $\h$-adic algebras $\alpha:\Uhg \rightarrow
  \Ug[[\h]]$ such that $\alpha = \id + \Ocal(\h)$ and $\alpha
  \rvert_{\mathfrak{h}} = \id_{\mathfrak{h}}$.
\end{Theorem}
\noindent
The important consequence of this theorem for representation theory is that weight vectors and weight spaces of representations of the deformed and undeformed algebras can be identified. While in this sense, the irreducible representations of quantum alebras are equivalent to the usual representations, the nonequivalent coproducts on the enveloping algebra and its quantum deformation lead to different tensor representation in the deformed and undeformed case. For Drinfeld-Jimbo deformations $\Uhsu$, which are the Hopf duals of quantum groups, the reduction of tensor representations are given by $q$-deformed Clebsch-Gordan coefficients. As the deformed and undeformed coproducts are related by a Drinfeld twist, it was quickly realized \cite{Curtright:1991} that the representations of Drinfeld twists ought to be given by a combination of deformed and undeformed Clebsch-Gordan coefficients. Indeed, one can rigorously proove the following Proposition \cite{Blohmann:2002a}:

\begin{Proposition}
  \label{th:Frep1}
  There is a universal Drinfeld twist $\F$ from $\Usu$ to $\Uhsu$, which has  the
  matrix elements
  \begin{equation}
  \label{eq:Frep}  
    \Braket{j_1,m_1';j_2,m_2'}{\F}{j_1,m_1;j_2,m_2}
    = \sum_{j,m}
    \CGq{j_1}{j_2}{j}{m_1'}{m_2'}{m} \,
    \CG{j_1}{j_2}{j}{m_1}{m_2}{m} 
  \end{equation}
  in an irreducible representation of $\Usu \otimes \Usu$ with weights $j_1$, $j_2$ and basis $\Ket{j_1,m_1;j_2,m_2} := \Ket{j_1,m_1} \otimes \Ket{j_2,m_2}$, where the expressions in parentheses denote the $q$-deformed and undeformed Clebsch-Gordan coefficients.
\end{Proposition}

Recall that the action of an enveloping algebra $\Ug$ on an algebra $\Xcal$ is called covariant if for all $x,y \in \Xcal$ and $g \in \Ug$
\begin{equation}
\label{eq:Covar1}
  g \tr (xy) = (g_{(1)} \tr x)(g_{(2)} \tr y) \,.
\end{equation}
In mathematical terminology $\Xcal$ is called a module algebra. For the undeformed coproduct Eq.~\eqref{eq:Covar2} simply means that the elements of the Lie algebra $\mathfrak{g} \subset \Ug$ act as derivations on $\Xcal$. A quantum space, which is by definition a module algebra of the quantum deformation $\Uhg$, is realized by a star product on a function algebra in a covariant manner only if the analogous condition 
\begin{equation}
\label{eq:Covar2}
  g \tr (x\star y) = (g_{(1_\h)} \tr x)\star (g_{(2_\h)} \tr y)
\end{equation}
holds, where $g_{(1_\h)} \otimes g_{(2_\h)} \equiv \Delta_\h(g)$ is the Sweedler notation for the deformed coproduct. If we define the star product map by
\begin{equation}
\label{eq:StarProd2}
  x \star y := (\F_{[1]} \tr x)(\F_{[2]} \tr y)\,,
\end{equation}
where we use the Sweedler like notation $\F_{[1]} \otimes \F_{[2]} \equiv \F$, covariance condition~\eqref{eq:Covar2} is satisfied because $\Delta_\h(g) = \F^{-1} \Delta(g) \F$. But are there twists for which Eq.~\eqref{eq:StarProd2} also defines an associative product, thus realizing the algebra of a quantum space? It turns out that the twist of Proposition~\ref{th:Frep1} realizes the quantum plane and, essentially, also quantum Euclidean 4-space, and quantum Minkowski space \cite{Blohmann:2002a}.

\section{The reconstruction method}
\label{sec:Method}

\subsection{The tensor operator basis}

We want to find a method to reconstruct elements of the enveloping algebra $\Usu$ from their irreducible representations. Consider the Cartan-Weyl basis $\{E,H,F\}$ of $\mathrm{su}_2$ with commutation relations
\begin{equation}
\label{eq:suCommutRel}
  [H,E]=2E\,,\qquad [H,F]=-2F \,,\qquad [E,F] = H \,,
\end{equation} 
the compact real form being given by the $*$-structure $E^* = F$, $H^* = H$, $F^* = E$. For our purposes the usual Poincar\'e-Birkhoff-Witt basis of ordered monomials of the generators
\begin{equation}
\label{eq:PBWBasis}
  \mathcal{B}_{\mathrm{PBW}}
  = \{E^i H^j F^k \,|\, i,j,k \in \mathbb{N}_0\}
\end{equation} 
turns out to be not particular convenient: The irreducible representations of the ordered monomials do not satisfy any obvious orthogonality relations which would allow us to draw immediate conclusions from the representations of a given algebra element to its coefficients with respect to this basis. Recall that for each half-integer weight $j \in \tfrac{1}{2} \mathbb{N}_0$ there is one irreducible unitary representation of $\Usu$ defined on the orthonormal weight-$j$ (or spin-$j$) basis $\{\Ket{j,m}, m= -j,-j+1,\ldots,j\}$ by
\begin{equation}
\label{eq:Irrepsu}
\begin{aligned}
  E\Ket{j,m} &= \sqrt{(j+m+1)(j-m)} \,\Ket{j,m+1} \\
  F\Ket{j,m} &= \sqrt{(j+m)(j-m+1)}\,\Ket{j,m-1}\\
  H\Ket{j,m} &= 2m\Ket{j,m} \,,
\end{aligned}
\end{equation}
The structure homomorphism $\rho^j: \Usu \rightarrow \mathrm{End}(\mathbb{C}^{2j+1})$ is given by the matrix elements, $\rho^j(g)^{m'}{}_{m} := \Braket{j,m'}{g}{j,m}$. Since the Lie algebra $\mathrm{su}_2$ is simple, any representation of $\Usu$ can be decomposed into a direct sum of irreducible subrepresentations, each of which is isomorphic to a representation given by~\eqref{eq:Irrepsu}. This is in particular true for the adjoint action of $\Usu$ on itself which is defined on the generators as 
\begin{equation}
\label{eq:adDef}
  \ad g \tr u := [g,u]\,,\qquad
  g\in \mathrm{su}_2 \subset \Usu \,,\quad u \in \Usu \,.
\end{equation}
Let $\{ T^j_m \in \Usu \,|\, m= -j,\ldots,j \}$ be a weight basis of a weight-$j$ subrepresentation of the adjoint representation, that is,
\begin{equation}
\label{eq:TensDef}
  [g, T^j_m] = \sum_{m'} T^J_{m'} \rho^j(g)^{m'}{}_{m}
  \equiv \sum_{m'} T^J_{m'} \Braket{j,m'}{g}{j,m}
\end{equation} 
for all $g \in \mathrm{su}_2$. Such a basis $\{T^j_m\}$ is called a weight-$j$ tensor operator of $\mathrm{su}_2$. The set of all weight-$0$ operators is the center of $\Usu$. As commutative algebra, the center is generated by the canonical quadratic Casimir element $C := \sum_{ij} g_i g_j K^{ij}$, where $\{g_i\}$ is a basis of the Lie algebra and $K^{ij}$ is the inverse of the Killing metric $K_{ij} := \mathrm{tr}(\ad g_i \,\ad g_j)$. In  the Cartan-Weyl basis we obtain
\begin{equation}
  C = \tfrac{1}{2}EF+ \tfrac{1}{2}FE+\tfrac{1}{8}H^2 
  = EF + \tfrac{1}{8}H(H-2) \,,
\end{equation}
such that the polynomial algebra $\mathbb{C}[C]$ is the center of $\Usu$. The representations of the Casimir element,
\begin{equation}
\label{eq:Casirep}
  C\Ket{j,m} = j(j+1)\Ket{j,m} \,,
\end{equation}
show that $C$ is the ususal square of angular momentum. 

By definition, $T^j_j$ is the highest weight vector of a weight-$j$ subrepresentation of the adjoint representation, so $(\ad E) \tr T^j_j \equiv [E,T^j_j]= 0$ and $(\ad C) \tr T^j_j = j(j+1) T_j^j$. From these two equalities it follows, that $T^j_j = z E^j$, where $z$ is some element of the center. If we pick $z$ from the number field we get the tensor operators
\begin{equation}
\label{eq:minimalTensor}
  T^J_J := \alpha E^J \,,\qquad \alpha \in \mathbb{C} \,,
\end{equation}
from which all other tensor operators can be obtained by multiplication by a central element. Here $\alpha$ is a normalization constant, which will later be chosen for convenience. From now on we denote by $T^J_M$ always the tensor operator which is generated by $\alpha E^J$. We use capital letters for the indices in order to allow in the formulas which we will derive below for a clear disctinction of the weights pertaining to the adjoint action from those pertaining to matrix representations. The fact that, as module with respect to the adjoint action, $\Usu$ can be completely decomposed into irreducible submodules implies that
\begin{equation}
\label{eq:TensorBasis}
  \mathcal{B}_{\mathrm{tensor}}
  = \{C^k T^J_M \,|\, k,J \in \mathbb{N}_0 \,, M = -J,-J+1,\ldots,J \}
\end{equation}
is a basis of $\Usu$, which we will call the tensor basis. The fact that~\eqref{eq:TensorBasis} is a basis of $\Usu$ means that the tensor operators are a basis of $\Usu$ as free module over its center. Thus, every element $a\in \Usu$ can be written uniquely as
\begin{equation}
\label{eq:Reconst1}
  a = \sum_{J,M} a^J_M\, T^J_M\,,  
  \qquad a^J_M \in \mathbb{C}[C] \,,
\end{equation}
where the sum runs over a finite subset of all allowed integer values of $J$ and $M$. Reconstructing the element $a$ from its representations now amounts to finding the polynomials $a^J_M$.

\subsection{The reconstruction method}

Let us compute the irreducible representations of Eq.~\eqref{eq:Reconst1}. First, we consider the central coefficients $a^J_M$. Since $a^J_M$ is a polynomial in the Casimir, the matrix element is a polynomial of the weight $j$ of the representation,
\begin{equation}
\label{eq:Reconst2}
  \Braket{j,m}{a^J_M}{j,m} =: a^J_M(j) \in \mathbb{C}[j]\,.
\end{equation}
Due to Eq.~\eqref{eq:Casirep} this polynomial satisfies
\begin{equation}
\label{eq:Reconst3}  
  a^J_M(j) = a^J_M(-j-1) \,,
\end{equation}
since it is actually a polynomial in $j(j+1)$ or, equivalently, a quadratic polynomial in $j+\tfrac{1}{2}$. Conversely, given a polynomial $p(j) \in \mathbb{C}[j]$ which satisfies $p(j) = p(-j-1)$ there exists a unique polynomial in the Casimir which has $p(j)$ as its matrix elements. For an intuitive notation we will denote this polynomial by $p(\jhat) \in \mathbb{C}[C]$, such that its defining equation takes the suggestive form
\begin{equation}
\label{eq:Reconst4}  
  \Braket{j,m}{p(\jhat)}{j,m} = p(j) \,.
\end{equation}
The map $p(j) \mapsto p(\jhat)$ could be viewed as substituion
\begin{equation} 
  j \mapsto \jhat = \tfrac{1}{2}(\sqrt{4C+1} -1) \,,
\end{equation}
where the relation $p(j) = p(-j-1)$ guarantees that the square roots drop out such that $p(\jhat)$ is a polynomial in $C$ only. We emphasize that we do not add such an square root of the Casimir to the algebra, though. We view $p(\jhat)$ merely as a suggestive notation for the element of the center which is uniquely defined by Eq.~\eqref{eq:Reconst4}.  

The matrix elements of the tensor operators are given by the Wigner-Eckhart theorem,
\begin{equation}
\label{eq:WignerEckart}
  \Braket{j,m'}{T^J_M}{j,m}
  = \rBraket{j}{T^J}{j} \CG{J}{j}{j}{M}{m}{m'} \,,
\end{equation}
where the reduced matrix element $\rBraket{j}{T^J}{j}$ does not
depend on $m$, $m'$, or $M$, and where the expression in parentheses denotes the Clebsch-Gordan coeffient. The explicit form and some properties of Clebsch-Gordan coefficients and their $q$-deformations can be found for example in \cite{Schmuedgen}. The reduced matrix elements will be computed below.

The irreducible representations of Eq.~\eqref{eq:Reconst1} now take the form
\begin{equation}
\label{eq:Reconst4b}
  \Braket{j,m'}{a}{j,m} = \sum_{J,M} a^J_M(j)\,
  \rBraket{j}{T^J}{j} \CG{J}{j}{j}{M}{m}{m'} \,.
\end{equation}
The main advantage of using the tensor basis~\eqref{eq:TensorBasis} instead of the Poincar\'e-Birkhoff-Witt basis~\eqref{eq:PBWBasis} is the fact, that the Clebsch-Gordan coefficients satisfy orthogonality relations which can be used in order to solve Eq.~\eqref{eq:Reconst4b} for $a^J_M(j)$. Using the well-known orthogonality relation
\begin{equation}
  \sum_{m,m'}\CG{J}{j}{j}{M}{m}{m'} \CG{J'}{j}{j}{M'}{m}{m'}
  = \frac{2j+1}{2J+1}\,\delta_{JJ'}\delta_{MM'}
\end{equation}
we thus arrive at 

\begin{Proposition}
\label{th:Reconst1}
Let $a \in \Usu$ be an element of the enveloping algebra with matrix elements $\Braket{j,m'}{a}{j,m}$ with respect to the irreducible representations defined in Eqs.~\eqref{eq:Irrepsu}. Let $T^J_M \in \Usu$ be the minimal degree tensor operators generated by $T^J_J \sim E^J$ and $\rBraket{j}{T^J}{j}$ their reduced matrix elements. Then
\begin{itemize}

\item[(i)]
For all integers $J \geq 0$ and $M$, $|M| \leq J$ the expression
\begin{equation}
\label{eq:Reconst5}
  a^J_M(j)
  :=
  \frac{(2J+1)}{ (2j+1)\rBraket{j}{T^J}{j} }
  \sum_{m,m'}\CG{J}{j}{j}{M}{m}{m'}
  \Braket{j,m'}{a}{j,m}
\end{equation}
defines  a polynomial in $j$ which is nonzero for only a finite number of values of $J$ and $M$. 

\item[(ii)]
The polynomials $a^J_M(j)$ are quadratic in $j+\tfrac{1}{2}$, the substitution
\begin{equation}
  (j_1+\tfrac{1}{2})^2 \mapsto C + \tfrac{1}{4}
\end{equation}
yielding polynomials in the Casimir element $C$ which are denoted by $a^J_M(\jhat)$.

\item[(iii)] The element $a$ can be written as
\begin{equation}
\label{eq:Reconst6}
  a = \sum_{J,M} a^J_M(\jhat)\, T^J_M\,.
\end{equation}
\end{itemize}

\end{Proposition}

This reconstruction method can be readily generalized to the tensor product $\Usu \otimes \Usu$: Let $a \in \Usu \otimes \Usu$ be an element of the tensor product, let $\Braket{j_1,m_1';j_2,m_2'}{a}{j_1,m_1;j_2,m_2}$ denote its matrix elements with respect to irreducibles representation of each tensor factor. First we need to calculate
\begin{multline}
\label{eq:Reconst7}
  a^{J_1 J_2}_{M_1 M_2}(j_1,j_2)
  :=
  \frac{(2J_1+1)(2J_2+1)}{(2j_1+1)(2j_2+1)
  \rBraket{j_1}{T^{J_1}}{j_1}\rBraket{j_2}{T^{J_2}}{j_2} } \\
  \sum_{\substack{m_1,m'_1\\m_2,m'_2}}
  \CG{J_1}{j_1}{j_1}{M_1}{m_1}{m'_1}\CG{J_2}{j_2}{j_2}{M_2}{m_2}{m'_2} 
  \Braket{j_1,m_1';j_2,m_2'}{a}{j_1,m_1;j_2,m_2} \,,
\end{multline}
which defines polynomials, which are quadratic in $(j_1+\frac{1}{2})$ and $(j_2+\frac{1}{2})$. Then we substitute
\begin{equation}
\label{eq:Reconst7a}
  (j_1+\tfrac{1}{2})^2 \mapsto (C + \tfrac{1}{4})\otimes 1 
  \,,\qquad
  (j_2+\tfrac{1}{2})^2 \mapsto 1 \otimes (C + \tfrac{1}{4})
\end{equation} 
in order to obtain the unique central elements   
\begin{equation}
  a^{J_1 J_2}_{M_1 M_2}(\jhat_1,\jhat_2)
  \in  \mathbb{C}[ C\otimes 1, 1\otimes C ] \,,
\end{equation}
the representations of which are given by the polynomials~\eqref{eq:Reconst7}. Finally, reconstruct the element of the tensor algebra by 
\begin{equation}
\label{eq:Reconst7b}
  a = \sum_{\substack{J_1,M_1\\ J_2, M_2}} 
  a^{J_1 J_2}_{M_1 M_2}(\jhat_1,\jhat_2)\,\,
  T^{J_1}_{M_1} \otimes T^{J_2}_{M_2} \,.
\end{equation}
We will now apply this reconstruction method to the Drinfeld twist~\eqref{eq:Frep}.

\section{Reconstruction of the universal Drinfeld twist}
\label{sec:Reconstruction}

\subsection{Calculation of the tensor basis}

In order to obtain explicit formulas from the reconstruction method 
of Proposition~\ref{th:Reconst1} we need to calcualte the tensor operators $T^J_M$ in terms of the Poincar\'e-Birkhoff-Witt basis as well as the reduced matrix elements $\rBraket{j}{T^J}{j}$. We start with the reduced matrix elements. 

From Eq.~\eqref{eq:Irrepsu} we can derive the representation of powers of the generators
\begin{equation}
\label{eq:EJrep}
\begin{aligned}
  E^J\Ket{j,m} &=
  \sqrt{(-1)^J (j+m+1)_J (m-j)_J} \,\Ket{j,m+J} \\
  F^J\Ket{j,m} &=
  \sqrt{(-1)^J (j-m+1)_J (-m-j)_J} \,\Ket{j,m-J}  \,,
\end{aligned}
\end{equation}
where
\begin{equation}
  (x)_J := (x)(x+1)\cdots(x+J-1)
\end{equation} 
denotes the Pochhammer symbol. From Eqs.~\eqref{eq:EJrep} we obtain for the irreducible representations of the tensor operator~\eqref{eq:minimalTensor} on the one hand
\begin{equation}
\label{eq:Red1}
  \Braket{j,j}{T^J_J}{j,j-J} = 
  \Braket{j,j}{\alpha E^J}{j,j-J} =  
  \alpha\sqrt{\frac{(2j)!\, J!}{(2j-J)!}} 
\end{equation}
for $J\leq 2j$. On the other hand we have due to the Wigner-Eckhart theorem~\eqref{eq:WignerEckart}
\begin{equation}
\label{eq:Red2}
\begin{split}
  \Braket{j,j}{T^J_J}{j,j-J}
  &= \rBraket{j}{T^J}{j} \CG{J}{j}{j}{J\,}{j-J}{j} \\
  &= \rBraket{j}{T^J}{j} 
  \sqrt{\frac{(2j+1)!(2J)!}{(2j+J+1)!\,J!}} \,\,\,,
\end{split}
\end{equation}
where we have inserted the explicit expression for the Clebsch-Gordan coefficient. We conclude that
\begin{equation}
\label{eq:Red3}
  \rBraket{j}{T^J}{j} =  \alpha
  \sqrt{ \frac{(2j+J+1)!\,J!\,J!}{(2j+1)(2j-J)!(2J)!} } \,\,\,.
\end{equation}
For our purposes, it is convenient to chose the normalization constant $\alpha$ such that
\begin{equation}
\label{eq:Red4}
  \rBraket{J}{T^J}{J} = 1 \,,
\end{equation}
for which we have to set
\begin{equation}
\label{eq:Red5}
  \alpha :=
  \sqrt{ \frac{(2J+1)!}{(3J+1)!\,J!} } \,\,\,.
\end{equation}
From now on we will assume this choice of $\alpha$, for which the value of the reduced matrix element~\eqref{eq:Red3} becomes
\begin{equation}
\label{eq:Red6}
  \rBraket{j}{T^J}{j} 
  = \sqrt{ \frac{(2J+1)(2j+J+1)!\,J!}{(2j+1)(3J+1)!(2j-J)!} } \,\,\,.
\end{equation}
From the heighest weight vector $T^J_J$ we obtain the weight basis by repeated action of the lowering operator $\ad F$. More precisely, from Eq.~\eqref{eq:EJrep} we conclude that
\begin{equation}
\label{eq:Red7}
\begin{split}
  T^J_{M} 
  &=  [(-1)^{J-M} (J-M)! (-2J)_{J-M}]^{-\frac{1}{2}}
  \,(\ad F)^{J-M} \tr T^J_J \\ 
  &= \sqrt{\frac{(2J+1)(J+M)!}{(3J+1)!\,J!\,(J-M)!}} 
  \,(\ad F)^{J-M} \tr E^J 
\end{split}
\end{equation}
for $|M| \leq J$. The remaining computational problem for an explicit expression of $T^j_m$ in terms of the Poincar\'e-Birkhoff-Witt basis is the lexicographic reordering of $(\ad F)^{J-M} \tr E^J$. Details of the computation are provided in Appendix~\ref{sec:AppTensor}. As result we obtain
\begin{subequations}
\label{eq:Red8}
\begin{align}
  T^J_M
  &= (-1)^{J+M} \sqrt{\frac{(2J+1)J!(J-M)!(J+M)!}{(3J+1)!}} \notag \\ 
  &\quad\times
  \sum_{p = 0}^{p\leq \frac{J-M}{2} } 
  \frac{(-1)^p}{p! (p+M)!}
  E^{p+M}\binom{J+H-1}{J-M-2p} F^p \qquad\text{for}\quad M \geq 0  
  \label{eq:Red8a}\\
  T^J_M
  &= (-1)^{J-M} \sqrt{\frac{(2J+1)J!(J-M)!(J+M)!}{(3J+1)!}} \notag \\ 
  &\quad\times
  \sum_{p = 0}^{p\leq \frac{J+M}{2} } 
  \frac{(-1)^p}{p! (p-M)!}
  E^{p}\binom{J+H-1}{J+M-2p} F^{p-M} \qquad\text{for}\quad M < 0 \,,
  \label{eq:Red8b}
\end{align}
\end{subequations}
where the algebra valued binomial coefficient is defined by
\begin{equation}
  \binom{X}{k} := \frac{(-1)^k(-X)_k}{k!} \,,
\end{equation}
denoting a polynomial in $X$. For a complete expansion in terms of ordered monomials we yet have to expand the binomials in powers of $H$,
\begin{equation}
\label{eq:Red9}
  \binom{J+H-1}{J \pm M-2p}
  = \sum_{n=0}^{J \pm M-2p} H^n \sum_{k=n}^{J-M \pm 2p}
  \frac{1}{k!} \binom{J-1}{J \pm M-2p-k} s(k,n) \,,
\end{equation}
where $s(k,n)$ are Stirling numbers of the first kind.

\subsection{The generating function for the Drinfeld twist}

We will now apply the reconstruction method of Sec.~\eqref{sec:Method} to the universal Drinfeld twist of $\F$ of Proposition~\ref{th:Frep1}. Inserting the representations~\eqref{eq:Frep} of the twist into Eq.~\eqref{eq:Reconst7} the twist can be expressed according to Eq.~\eqref{eq:Reconst8} as
\begin{equation}
\label{eq:Reconst8}
\begin{split}
  \F =&
  \sum_{\substack{J_1,M_1\\ J_2, M_2}}
  \frac{(2J_1+1)(2J_2+1)}{(2\jhat_1+1)(2\jhat_2+1)
  \rBraket{\jhat_1}{T^{J_1}}{\jhat_1}\rBraket{\jhat_2}{T^{J_2}}{\jhat_2} } \\
  &\sum_{\substack{m_1,m'_1\\m_2,m'_2}}
    \CG{J_1}{\jhat_1}{\jhat_1}{M_1}{m_1}{m'_1}
    \CG{J_2}{\jhat_2}{\jhat_2}{M_2}{m_2}{m'_2} \\ 
  &\sum_{j,m}
    \CG{\jhat_1}{\jhat_2}{j}{m'_1}{m'_2}{m}
    \CGq{\jhat_1}{\jhat_2}{j}{m_1}{m_2}{m}
  \,\,T^{J_1}_{M_1} \otimes T^{J_2}_{M_2}
\end{split}
\end{equation}
where we recall that the hats on $\jhat_1$ and $\jhat_2$ indicate that the coefficients of the tensor operators are polynomials in $C\otimes 1$ and $1 \otimes C$ which we obtain after substitution~\eqref{eq:Reconst7a}. Eq.~\eqref{eq:Reconst8} does in general not yield an element of $\Usu \otimes \Usu$ for any fixed value of $q$. It has to be understood as algebra valued generating function in $q=\E^\h$ which produces in each order of $\h$ an element of $\Usu \otimes \Usu$ proper. An explicit expansion up to third order in $\h$ will be given in the next section.

Note that while the entire dependence on $\h$ is contained in the $q$-deformed Clebsch-Gordan coefficient, the arguments of the latter are contracted with the arguments of undeformed Clebsch-Gordan coefficients. We can confine the $\h$-de\-pen\-dence further by using the following identity for the Clebsch-Gordan coefficients which is derived in Appendix~\ref{sec:AppNineJ}:  
\begin{multline}
\label{eq:Reconst9}  
  \CG{J_1}{j_1}{j_1}{M_1}{m_1}{m'_1}
  \CG{J_2}{j_2}{j_2}{M_2}{m_2}{m'_2} 
  \CG{j_1}{j_2}{j}{m'_1}{m'_2}{m}  
  =\\
  \sum_{J,j'}
  \beta
  \NineJ{J_1}{j_1}{j_1}  
        {J_2}{j_2}{j_2}
        {J  }{j' }{j  }
  \CG{J_1}{J_2}{J}{M_1}{M_2}{M}
  \CG{j_1}{j_2}{j'}{m_1}{m_2}{m'} 
  \CG{J}{j'}{j}{M}{m'}{m} \,,
\end{multline}
where the expression in braces denotes the $9j$-symbol, the factor $\beta$ is defined as
\begin{equation}
\label{eq:Reconst10}  
  \beta :=
  \sqrt{(2J+1)(2j'+1)(2j_{1}+1)(2j_{2}+1)} \,,
\end{equation}
and $m'_1 = M_1 + m_1$, $m'_2 = M_2 + m_1$, $M = M_1 + M_2$, $m' = m_1 + m_2$. Inserting Eq.~\eqref{eq:Reconst9} into Eq.~\eqref{eq:Reconst8} we obtain
\begin{equation}
\begin{split}
\label{eq:Reconst11}
  \F
  =&
  \sum_{J_1,J_2,J}\,\, \sum_{j,j'}
  \frac{(2J_1+1)(2J_2+1)}{
  \rBraket{\jhat_1}{T^{J_1}}{\jhat_1}
  \rBraket{\jhat_2}{T^{J_2}}{\jhat_2} } 
  \sqrt{\frac{(2J+1)(2j'+1)}{(2\jhat_1+1)(2\jhat_2+1)}}
  \NineJ{J_1}{\jhat_1}{\jhat_1}
        {J_2}{\jhat_2}{\jhat_2}
        {J  }{j' }{j  } 
  \\
  &\times \sum_m \CG{J}{j'}{j}{0}{m}{m}  
  \sum_{m_1,m_2}
    \CG{\jhat_1}{\jhat_2}{j'}{m_1}{m_2}{m}
    \CGq{\jhat_1}{\jhat_2}{j}{m_1}{m_2}{m}    
  \\
  &\times \sum_{M} \CG{J_1}{J_2}{J}{M}{-M}{0}
  \,\, T^{J_1}_{M} \otimes T^{J_2}_{-M} 
\end{split}
\end{equation}
where we have used that from condition $m'_1 + m'_2 = m= m_1 + m_2$ in Eq.~\eqref{eq:Reconst8} it follows that $M = 0$ and $m' = m$. 

In the form of Eq.~\eqref{eq:Reconst11} the generating function gives us some insight into the structure of the the twist. The first line of Eq.~\eqref{eq:Reconst11} and the first Clebsch-Gordan coefficient on the second line do not depend on the deformation parameter $\h$ and contain only well known functions, the $9j$-symbol and the Clebsch-Gordan coefficient essentially being given by hypergeometric functions. 

The summation over $M$ in last line eliminates the dependence on the magnetic quantum numbers $M_1$ and $M_2$ of the tenor operator basis. The fact that the magnetic quantum number of the tensor operators $T^{J_1}_{M_1} \otimes T^{J_2}_{M_2}$ which appear in the Drinfeld twist add up to zero, $M_1 + M_2 = M = 0$, can also be understood on a more abstract level: Up to isomorphism, the quantum deformation of an enveloping algebra does not affect the Cartan subalgebra as it was stated in Theorem~\ref{th:CartanPreserve}. For the Drinfeld-Jimbo deformation $\Uhsu$ which we consider here this means that
\begin{equation}
  \Delta(H) = \Delta_\hbar(H) = H\otimes 1 + 1 \otimes H \,,
\end{equation}
which implies that the Drinfeld twist $\F$ must commute with $\Delta(H)$. From
\begin{equation}
  [\Delta(H), T^{J_1}_{M_1} \otimes T^{J_2}_{M_2} ]
  = 2(M_1 + M_2)(T^{J_1}_{M_1} \otimes T^{J_2}_{M_2})
\end{equation} 
we conclude that only those products of tensor operators can appear in $\F$ for which $M_1 + M_2= 0$. 

The dependence of Eq.~\eqref{eq:Reconst11} on the deformation parameter is contained in  the contraction of the deformed and undeformed Clebsch-Gordan coefficient over $m_1$ and $m_2$ in the second line. The representation theoretic interpretation of this term is the following: We can use both, the undeformed and the deformed coproduct, to define a tensor product representation of two irreducible representations with weights $j_1$ and $j_2$, defining the undeformed and deformed strucure maps as
\begin{equation}
  \rho^{j_1 \otimes j_2} := 
  (\rho^{j_1} \otimes \rho^{j_2}) \circ \Delta
  \quad\text{and}\quad
  \rho_\h^{j_1 \otimes j_2} := 
  (\rho^{j_1} \otimes \rho^{j_2}) \circ \Delta_\h \,.
\end{equation}
Both representations can be reduced into irreducible components. Denoting the basis vectors of the irreducible weight-$j$ subrepresentation of the undeformed and deformed tensor representation by $\Ket{j_1, j_2 \rightarrow j,m}$ and $\Ket{j_1,j_2 \rightarrow j,m}_{\h}$, respecively, we obtain
\begin{equation}
\label{eq:hybrid}  
  \langle j_1,j_2 \rightarrow j',m \,\vert\,
  j_1, j_2 \rightarrow j,m \rangle_{\h}
  = \sum_{m_1,m_2}
    \CG{j_1}{j_2}{j'}{m_1}{m_2}{m}
    \CGq{j_1}{j_2}{j}{m_1}{m_2}{m} \,.
\end{equation}  
In other words the deformation is now expressed as the change of basis from the irreducible components of tensor representations with respespect to the undeformed coproduct $\Delta$ to those with respect to the deformed coproduct $\Delta_\h$. Again, the expression on the right hand side of Eq.~\eqref{eq:hybrid} is to be understood as generating function. While the $q$-Clebsch-Gordan coefficient are well known functions for a given value of $q$, little is known about its perturbative expansion in powers of $\h$.

\subsection{Perturbative expansion}
\label{sec:Perturb}

\newcommand{\qZahl}[1]{(#1)_q}
\newcommand{\qPoch}[2]{(#1;q)_{#2}}

Ideally, we would like to find a closed form expression for the Drinfeld twist in each order of $\h$. This would require a closed form expansion of Eq.~\eqref{eq:hybrid}, which is essentially given by a sum of the product of the ordinary hypergeometric function ${}_3 F_2$ and its basic ($q$-deformed) counterpart ${}_3 \varphi_2$. To our best knowledge such hybrid combinations of ordinary and basic hypergeometric functions have not been studied in the literature yet and little is known about the perturbative expansion of basic hypergeometric functions in powers of $\hbar = \ln q$ or other possible perturbation parameters such as $q - q^{-1}$ and $q-1$. Studying the general problem of perturbative expansion of basic hypergeometric functions is beyond the scope of this article. This is ongoing research and will be presented elsewhere. Here we will only expand the $q$-deformed Pochhammer symbol, which is the building block of basic hypergeometric functions. This will enable us to carry out the explicit caculation of each order of the Drinfeld twist by a Taylor series expansion of the generating functions~\eqref{eq:Reconst8} and~\eqref{eq:Reconst11}.

For our puroses it is convenient to consider the $q$-Pochhammer symbol $[x]_n$ which is defined by symmetric quantum numbers $[x]$,
\begin{equation}
\label{eq:Perturb1}  
  [x]_n := [x]\cdot [x+1]\cdots [x+n-1]
  \,,\qquad
  [x] = \frac{\E^{x\hbar}-\E^{-x\hbar}}{\E^{\hbar}-\E^{-\hbar}}
  =\frac{\sinh x\hbar}{\sinh \hbar} \,.
\end{equation}
Considering the logarithm of the Pochhammer symbols will turn the product of the $q$-numbers into the sum of their logarithms. Using the well known formula
\begin{equation}
\label{eq:Perturb2}  
  \ln \frac{\sinh x}{x}
  = \sum_{k=1}^{\infty} \frac{2^{2k-1} B_{2k}}{k(2k)!} 
  \,x^{2k} \,,
\end{equation}
where $B_{2k}$ are Bernoulli numbers, we obtain for the expansion of the logarithm of a quantum number
\begin{equation}
\label{eq:Perturb3}  
  \ln \frac{[x]}{x}
  = 
  \sum_{k=1}^{\infty} \frac{2^{2k-1} B_{2k}}{k(2k)!}
  \,(x^{2k}-1) \hbar^{2k} \,.
\end{equation} 
The Pochhammer symbol can then be expressed as exponential of the sum of this power series,
\begin{equation}
\label{eq:Perturb4}  
  \frac{[x]_n}{(x)_n}
  = \exp\left(
  \sum_{k=1}^{\infty} \frac{2^{2k-1} B_{2k}}{k(2k)!}
  \,\sum_{j=0}^{n-1}\bigl\{ (x+j)^{2k}-1 \bigr\} \hbar^{2k} \right) \,.
\end{equation} 
The sum is carried out using
\begin{equation}
\label{eq:Perturb5}  
  \sum_{j=0}^{n-1} (x+j)^{2k} = \frac{B_{2k+1}(x+n)-B_{2k+1}(x)}{2k+1} \,,
\end{equation}
where $B_{k}(x)$ denotes Bernoulli polynomials. We thus get
\begin{equation}
\label{eq:Perturb6}
  \frac{[x]_n}{(x)_n}
  = \exp\left(
  \sum_{k=1}^{\infty} \frac{2^{2k-1} B_{2k}}{k(2k)!}
  \biggl\{ \frac{B_{2k+1}(x+n)-B_{2k+1}(x)}{2k+1} - n \biggr\} 
  \hbar^{2k} \right) \,.
\end{equation}
This formula could serve as starting point for a perturbative expansion of general $q$-hypergeometric functions. Here it suffices to deduce from Eq.~\eqref{eq:Perturb6} the expansion of the Pochhammer symbol in $\h$. Up to third order we obtain
\begin{equation}
\label{eq:Perturb7}
\begin{split}
  \frac{[x]_n}{(x)_n}
  &= 1 + 
  \tfrac{1}{3}{B_{2}}\{B_{3}(x+n)-B_{3}(x) - 3 n \} \hbar^{2}
  + \Ocal(\h^4) \\
  &= 1 + \tfrac{1}{36}(-5\,n - 3\,n^2 + 2\,n^3 - 
    6\,n\,x + 6\,n^2\,x + 6\,n\,x^2) \h^2 + \Ocal(\h^4) \,.
\end{split}
\end{equation}
This expression is polynomial in $n$ and $x$ to each order of $\h$. Inserting it into 
Eq.~\eqref{eq:Reconst7} yields the searched-for polynomials in $C \otimes 1$ and $1 \otimes C$. From the generating functions~\eqref{eq:Reconst8} or~\eqref{eq:Reconst11} we then obtain the universal Drinfeld twist up to third order in $\h$. 

Each order $\F_k$ of the expansion $\F = \sum_k \F_k \h^k$ is alternatingly symmetric or antisymmetric with respect to the exchange of tensor factors by the transpose $\tau (a \otimes b) = b \otimes a$ according to
\begin{equation}
  \tau(\F_k) = (-1)^k \F_k \,. 
\end{equation}
This property can be derived from the fact that the transpose of the deformed coproduct amounts to a change of sign of the perturbation parameter, from which it follows that $ \tau( \F(\h)) = \F(-\h)$. Alternatively, it can be derived from the symmetry properties of the $q$-Clebsch-Gordan coefficients with respect to the transformation $q \mapsto q^{-1}$. The explicit expressions for the first three orders of the twist we finally obtain are
\begin{subequations}
\label{eq:Fresult}
\begin{align}
  \F_1 &= 
  2(T^{1}_{-1}\otimes T^{1}_{1} - T^{1}_{-1}\otimes T^{1}_{1})
  = 2\, T^{1}_{-1}\otimes T^{1}_{1} - \text{transpose} \\
  \F_2 &=
  -\frac{1}{18}\,{ C} \otimes { C} 
  +\frac{\sqrt{14}}{6}\,T^{2}_{0} \otimes { C} 
  + \frac{\sqrt{21}}{6}(
   T^{1}_{ 1}\otimes T^{2}_{-1}
  - T^{1}_{-1}\otimes T^{2}_{1}) \notag\\
  &+ \frac{21}{2} \, T^{2}_{-2}\otimes T^{2}_{2} 
  - \frac{7}{4}\, T^{2}_{0}\otimes T^{2}_{0} + \text{transpose} \\
  \F_3 &=
  \frac{\sqrt{2}}{180}\, (3-4C)T^{1}_{0} \otimes C 
  + \frac{\sqrt{7}}{30}\, T^{2}_{0} \otimes (9-2C)T^{1}_{0} \notag\\
  &+ \frac{1}{75}\bigl[
  7 - 21(C\otimes 1 + 1 \otimes C) - 12 \, C \otimes C \bigr]
  T^{1}_{-1}\otimes T^{1}_{1} 
  + \frac{7}{2}\, T^{2}_{-2}\otimes T^{2}_{2} \notag\\
  &+ \frac{\sqrt{6}}{3}\,T^{3}_{0}\otimes C 
  + \frac{2\sqrt{2}}{5} \bigl[
  T^{3}_{1} \otimes (1-3C)T^{1}_{-1} 
  - T^{3}_{-1} \otimes (1-3C)T^{1}_{1} \bigr] \notag \\
  &+ \sqrt{21}\, T^{2}_{0} \otimes T^{3}_{0}
  + \sqrt{105}( T^{3}_{-2} \otimes T^{2}_{2} 
  + T^{3}_{2} \otimes T^{2}_{-2}) \notag\\
  &+ 18 ( 5 T^{3}_{-3} \otimes T^{3}_{2} + T^{3}_{-1} \otimes T^{3}_{1})  
  - \text{transpose}
\end{align}
\end{subequations}
where ``transpose'' is shorthand for the tensor transpose of all preceeding terms such that each expression becomes symmetric or antisymmetric, respectively. One can use Eqs.~\eqref{eq:Red8} in order to express the result in terms of the Poincar\'e-Birkhoff-Witt basis. However, this yields expressions which are much longer than those of Eqs.~\eqref{eq:Fresult}, indicating that the tensor operator basis seems to be the better choice within the context of Drinfeld twists.  

The calculations leading to Eqs.~\eqref{eq:Fresult} are elementary but lengthy and are best carried out using computer algebra. With the expansion~\eqref{eq:Perturb6} of the $q$-Pochhammer symbol at hand the Taylor series expansion of the generating function~\eqref{eq:Reconst8} is reduced to addition and multiplication of polynomials, operations which are implemented efficiently by all common computer algebra systems. Hence, the explicit calculation of the Drinfeld twist to third order is not significantly limited by computing resources in any way. In any case, by the method presented here it is possible to compute the twist explicitly to orders which are high enough for the applications of Drinfeld twists to mathematical physics which we had in mind.

\section{Conclusion}
\label{sec:Conclusion}

Although the existence of universal Drinfeld twists can be proved rather easily, their calculation is a notoriously difficult and long standing problem. While we still did not derive a closed form for each order in the perturbation parameter of the universal twist of $\Uhsu$, significant progress towards this goal was presented here: In Eq.~\eqref{eq:Reconst11} we have given a generating function for the twist to all orders which can be easily expanded in powers of $\h$, as demonstrated in Eq.~\eqref{eq:Fresult}. Moreover, the generating function, which is expressed in terms of basic and ordinary hypergeometric functions, gives new insight into the general structure of the twist.

It is not difficult to understand why the proof of existence of the twist is so easy but the computation is so hard: The existence proof relies mainly on the fact that the first Hochschild cohomolgy of the enveloping algebra is zero. This means that every 1-cocycle is the coboundary of a 0-cocycle or, in other words, every derivation is inner. But we do not know how to compute this 0-cocycle. If in analogy to differential forms we view the inversion of the coboundary operator as a sort of integration, then the non-constructive existence proof uses integrability but does not tell us how to actually integrate. Just as in differential calculus, this cohomolgical type of integration turns out to be a difficult problem. In contrast, the series expansion of the generating function~\eqref{eq:Reconst11} in powers of $\h$ is a problem of differentiation. While integration is an art, differentiation is a simple technique which can be left to a computer algebra system. This is the reason why we consider the availability of a generating function as significant progress.

The computer algebra expansion of the generating function is computationally cheap and produces expansions of the twist which will suffice for many applications. However, it is not completely satisfactory as it produces expansion formulas like Eq.~\eqref{eq:Fresult} containing a lot of ``magical'' combinatorial numbers which cannot be explained any further. In Eq.~\eqref{eq:Reconst11} the dependence of the twist on the perturbation parameter is entirely confined to the $q$-Clebsch-Gordon coeffient, that is, essentially to the basic hypergeometric function ${}_3\varphi_2$ with basis $q=\E^\h$. Hence, the remaining problem which still separates us from a truly closed form expression for the universal twist is the perturbative expansion of this basic hypergeometric function in powers of $\h$. To our best knowledge, the question of perturbative expansion of basic hypergeometric functions, which seems so obvious in the context of quantum groups, has so far not recieved any systematic treatment in the special functions literature. Therefore, we had to make in Sec.~\ref{sec:Perturb} our own first step in this direction, computing a closed form expression for the $q$-deformed Pochhammer symbol in Eq.~\eqref{eq:Perturb6}.  We believe that, further pursuing this approach, a closed form expansion of basic hypergeometric functions and, hence, a closed form of the universal Drinfeld twist of $\Uhsu$ can be achieved.

\section*{Acknowledgements}

I would like to thank Daniel Sternheimer for hinting me to previous work on this subject. To Petr Kulish I am grateful for an illuminating discussion on the calculation of universal twists in the triangular case. I am also indebted to Hjalmar Rosengren for discussing aspects of the present work from the viewpoint of the theory of basic hypergeometric functions. 

\appendix

\section{Calculation of the tensor operator basis}
\label{sec:AppTensor}

Expressing the tensor operator basis in terms of the Poincar\'e-Birkhoff-Witt basis amounts to the normal ordering of Eq.~\eqref{eq:Red7}. While it is possible to carry out the normal ordering using the commutation relations of $\Usu$, this turns out to be surpisingly cumbersome. Therefore, we present an alternative approach which is much more in the spirit of this article: We deduce the normal ordered expression from the representations of the tensor operators.

Let us assume that $M\geq 0$. Starting from the Wigner-Eckart theorem~\eqref{eq:WignerEckart}, using~\eqref{eq:Red6} for the reduced matrix elements and the well-known explicit formula 
\begin{equation}
\begin{split}
  \CG{j_1}{j_2}{j}{m_1}{m_2}{m}
  &= (-1)^{m_1-j_1} 
  \sqrt{\frac{(2j+1)(j_1+j_2-j)!}{(j_1+j_2+j+1)!(j_1-j_2+j)!(j_2-j_1+j)!}}\\
  &\qquad \frac{(j_2+j-m_1)!}{(j_2-j+m_1)!}
  \sqrt{\frac{(j_1+m_1)!(j_2-m_2)!(j+m)!}{(j_1-m_1)!(j_2+m_2)!(j-m)!}}\\
  &\qquad
  {}_3 F_2\binom{m_1-j_1,\, j_1 + m_1 + 1,\, m-j}
  {j_2-j+m_1+1,\, -j-j_2+m_1 }
\end{split}
\end{equation}
for the Clebsch-Gordan coefficients \cite{Schmuedgen}, we derive for the matrix elements of the tensor operators
\begin{equation}
\label{eq:Tens1}
\begin{split}
  \Braket{j,m'}{T^J_M}{j,m}
  &= (-1)^{J+M} \sqrt{\frac{(2J+1)(J-M)!}{(3J+1)!(J)!(J+M)!}} \\
  &\quad \delta_{m',m+M} \sqrt{(-1)^M (j+m +1)_M (-j+m)_M}\\  
  &\quad  \sum_k (-1)^{k}\binom{J+M}{k} 
  (-j-m-k)_J (j+1-m-k)_J
\end{split}
\end{equation}
We want to deduce the element of the algebra in the Poincar\'e-Birkhoff-Witt basis from these representations. Towards this end we will compare Eq.~\eqref{eq:Tens1} with the matrix elements of monomials 
\begin{equation}
\begin{aligned}
\label{eq:Tens2}  
  \Braket{j,m'}{E^pF^p}{j,m} &= \delta_{m'm}(-1)^p(-j-m)_p(j-m+1)_p \\
  \Braket{j,m'}{E^M}{j,m} &= \delta_{m',m+M} \sqrt{(-1)^M(j+m+1)_M(-j+m)_M}
\end{aligned}
\end{equation}
We immediately see that the second line of Eq.~\eqref{eq:Tens1} is the matrix element of $E^M$. The last line has yet to be written in a different form. For this, we need a variant of the Pfaff-Saalsch\"utz summation formula
\begin{equation}
\label{eq:Saalschuetz}
  (a-c)_n (b-c)_n 
  = \sum_{p = 0}^n \binom{n}{p}
  (-c)_{n-p} (a+b-c+p)_{n-p}
  \, (a)_p (b)_p \,,
\end{equation} 
from which we get for $a = -j-m$, $b= j-m+1$, $c=k$, $n=J$
\begin{multline}
\label{eq:Summation1}
  (-j-m-k)_J (j-m+1-k)_J \\ 
  = \sum_{p = 0}^J \binom{J}{p}
  (-k)_{J-p} (-2m+1-k+p)_{J-p}
  \, (-j-m)_p (j-m+1)_p
\end{multline} 
and a variant of the Vandermonde summation formula
\begin{equation}
\label{eq:Summation2}
  \sum_{k=0}^{n}(-1)^k \binom{n}{k}
  (-k)_q (-k-x)_q\\
  = (-1)^n n! \binom{q}{n-q} (-x-q)_{2q-n} \,.
\end{equation}
Inserting first Eq.~\eqref{eq:Summation1} and then Eq.~\eqref{eq:Summation2} into the last line of Eq.~\eqref{eq:Tens1} we obtain
\begin{multline}
  \sum_{k=0}^{J+M} (-1)^{k}\binom{J+M}{k} 
  (-j-m-k)_{J} (j+1-m-k)_{J}\\
  = \sum_{p = 0}^{p\leq \frac{J-M}{2} } 
  \frac{J!(J+M)!}{p! (M+p)!}
  \binom{2m-2p +J - 1}{J-M-2p} (-j-m)_p (j-m+1)_p 
\end{multline}
Comparing this with the matrix element~\eqref{eq:Tens2}, we obtain the equality of matrix elements
\begin{multline}
  \Braket{j,m'}{T^J_M}{j,m}
   = (-1)^{J+M} \sqrt{\frac{(2J+1)J!(J-M)!(J+M)!}{(3J+1)!}}\\ 
  \sum_{p = 0}^{p\leq \frac{J-M}{2} } 
   \frac{(-1)^p}{p! (M+p)!} \,\,
   \Braket{j,m'}{E^{p+M}\binom{J+H-1}{J-M-2p} F^p}{j,m}
\end{multline}
from which we can deduce Eq.~\eqref{eq:Red8a}. The analogous calculations for $M \leq 0$ lead to Eq.~\eqref{eq:Red8b}.

\section{Derivation of Eq.~(\ref{eq:Reconst9})}
\label{sec:AppNineJ}

In order to derive Eq.~\eqref{eq:Reconst9} we recall that, while the Clebsch-Gordan coefficients reduce tensor representations, this reduction is neither commutative nor associative. Let us denote by $D^j$ the irreducible weight-$j$ representation. The isomorphism which corresponds to the exchange of the order in a product representation, $D^{j_1} \otimes D^{j_2} \rightarrow D^{j_1} \otimes D^{j_2}$, is given by a change of sign
\begin{equation}
\label{eq:NineJ1}
  \CG{j_1}{j_2}{j}
     {m_1}{m_2}{m}
  = (-1)^{j-j_1-j_2}
  \CG{j_2}{j_1}{j}
     {m_2}{m_1}{m} \,,
\end{equation}
where $j_1+j_2-j$ is always an integer. The associator which corresponds to changing the order of reduction of a product of three irreducible representations, $D^{j_1} \otimes (D^{j_2} \otimes D^{j_3})_{j_{23}} \rightarrow (D^{j_1} \otimes D^{j_2})_{j_{12}} \otimes D^{j_2}$ is by definition given by the Racah-coefficients
\begin{equation}
\label{eq:NineJ2}  
  \CG{j_2}{j_3}{j_{23}}
     {m_2}{m_3}{m_{23}}
  \CG{j_1}{j_{23}}{j}
     {m_1}{m_{23}}{m}
  = \sum_{j_{12}}
  \CG{j_1}{j_2}{j_{12}}
     {m_1}{m_2}{m_{12}}
  \CG{j_{12}}{j_{3}}{j}
     {m_{12}}{m_{3}}{m}
  \RCs{j_1}{j_2}{j_3}
      {j_{12}}{j_{23}}{j} \,,
\end{equation} 
where $m_{12}=m_1+m_2$, $m_{23}=m_2+m_3$. Using Eqs.~\eqref{eq:NineJ1} and \eqref{eq:NineJ2} the change of the reduction of a tensor product of four representations according to 
\begin{equation}
\label{eq:NineJ3}  
\begin{split}
  &(D^{j_1} \otimes D^{j_2})_{j_{12}} 
  \otimes (D^{j_1} \otimes D^{j_1})_{j_{34}} \\
  \rightarrow&
  ((D^{j_1} \otimes D^{j_2})_{j_{12}} 
  \otimes D^{j_3})_{j'} \otimes D^{j_4} \\
  \rightarrow&
  (D^{j_3} \otimes (D^{j_1} \otimes D^{j_2})_{j_{12}})_{j'} 
  \otimes D^{j_4} \\
  \rightarrow&
  ((D^{j_3} \otimes D^{j_1})_{j_{13}} \otimes D^{j_2})_{j'} 
  \otimes D^{j_4} \\
  \rightarrow&
  (D^{j_1} \otimes D^{j_3})_{j_{13}} \otimes 
  (D^{j_2} \otimes D^{j_4})_{j_{24}}   
\end{split}
\end{equation} 
is then expressed as
\begin{equation}
\label{eq:NineJ4}  
\begin{split}
  & 
  \CG{j_1}{j_2}{j_{12}}{m_1}{m_2}{m_{12}}
  \CG{j_3}{j_4}{j_{34}}{m_3}{m_4}{m_{34}} 
  \CG{j_{12}}{j_{34}}{j}{m_{12}}{m_{34}}{m}  
  \\
  = \sum_{j'}&
  \CG{j_1}{j_2}{j_{12}}{m_1}{m_2}{m_{12}}
  \CG{j_{12}}{j_3}{j'}{m_{12}}{m_3}{m'} 
  \CG{j'}{j_{4}}{j}{m'}{m_{4}}{m}  
  \RCs{j_{12}}{j_3}{j_4}{j'}{j_{34}}{j}
  \\
  = \sum_{j'}&
  \CG{j_1}{j_2}{j_{12}}{m_1}{m_2}{m_{12}}
  \CG{j_3}{j_{12}}{j'}{m_3}{m_{12}}{m'} 
  \CG{j'}{j_{4}}{j}{m'}{m_{4}}{m}  
  \RCs{j_{12}}{j_3}{j_4}{j'}{j_{34}}{j}
  (-1)^{j'-j_{12}-j_3}
  \\  
  = \sum_{j_{13},j'}&
  \CG{j_3}{j_1}{j_{13}}{m_3}{m_1}{m_{13}}
  \CG{j_{13}}{j_{2}}{j'}{m_{13}}{m_{2}}{m'} 
  \CG{j'}{j_{4}}{j}{m'}{m_{4}}{m} \\
  \times &
  \RCs{j_{3}}{j_{1}}{j_{2}}{j_{13}}{j_{12}}{j'}  
  \RCs{j_{12}}{j_3}{j_4}{j'}{j_{34}}{j}
  (-1)^{j'-j_{12}-j_3}
  \\    
  = \sum_{j_{13},j_{24}}&
  \CG{j_1}{j_3}{j_{13}}{m_1}{m_3}{m_{13}}
  \CG{j_{2}}{j_{4}}{j_{24}}{m_{2}}{m_{4}}{m_{24}} 
  \CG{j_{13}}{j_{24}}{j}{m_{13}}{m_{24}}{m} \\
  \times & \sum_{j'}
  \RCs{j_{13}}{j_{2}}{j_{4}}{j'}{j_{24}}{j}    
  \RCs{j_{3}}{j_{1}}{j_{2}}{j_{13}}{j_{12}}{j'}  
  \RCs{j_{12}}{j_3}{j_4}{j'}{j_{34}}{j}
  (-1)^{j'-j_{12}+j_{13}-j_1-2j_3} \,,
\end{split}
\end{equation}
where $m_{ij}=m_i+m_j$ for $i,j \in \{1,2,3,4\}$, $i<j$. Next we express the Racah-coefficients in terms of $6j$-symbols  
\begin{equation}  
\label{eq:NineJ5}
    \RCs{j_{1} }{j_{2} }{j_{3} }
        {j_{12}}{j_{13}}{j     }
  = (-1)^{j_1+j_2+j_3+j}\sqrt{(2j_{12}+1)(2j_{13}+1)}
  \SixJR{j_{1} }{j_{2} }{j_{3} }
        {j_{12}}{j_{13}}{j     } \,.
\end{equation}
Using the symmetries of the $6j$-symbol and the definition of the $9j$-symbol we can rewrite the last line of Eq.~\eqref{eq:NineJ4} as
\begin{equation}
\label{eq:NineJ6}
\begin{split}
  &\sum_{j'}
  \RCs{j_{13}}{j_{2}}{j_{4}}
      {j'}{j_{24}}{j}    
  \RCs{j_{3}}{j_{1}}{j_{2}}
      {j_{13}}{j_{12}}{j'}  
  \RCs{j_{12}}{j_3}{j_4}
      {j'}{j_{34}}{j} 
  (-1)^{j'-j_{12}+j_{13}-j_1-2j_3}
  \\
  =&\alpha \sum_{j'} (-1)^{2j'}(2j'+1)
  \SixJR{j_{13}}{j_{2}}{j_{4}}{j'}{j_{24}}{j}    
  \SixJR{j_{3}}{j_{1}}{j_{2}}{j_{13}}{j_{12}}{j'}  
  \SixJR{j_{12}}{j_3}{j_4}{j'}{j_{34}}{j}\\
  =&\alpha \sum_{j'} (-1)^{2j'}(2j'+1)
  \SixJ{j_{3} }{j_{13}}{j_1   }
       {j_2   }{j_{12}}{j'    }    
  \SixJ{j_{4} }{j_{24}}{j_{2} }
       {j_{13}}{j'    }{j     }  
  \SixJ{j_{34}}{j     }{j_{12}}
       {j'    }{j_{3} }{j_4   } \\
  =&\alpha
  \NineJ{j_{1} }{j_{2} }{j_{12} }
        {j_{3} }{j_{4} }{j_{34}}
        {j_{13}}{j_{24}}{j     }	
\end{split}
\end{equation} 
where the factor $\alpha$ is defined as
\begin{equation}
\label{eq:NineJ7}  
  \alpha := (-1)^{2j} 
  \sqrt{(2j_{12}+1)(2j_{34}+1)(2j_{13}+1)(2j_{24}+1)} \,.
\end{equation}
From Eqs.~\eqref{eq:NineJ4} and \eqref{eq:NineJ6} we finally obtain
\begin{multline}
\label{eq:NineJ8}  
  \CG{j_1}{j_2}{j_{12}}{m_1}{m_2}{m_{12}}
  \CG{j_3}{j_4}{j_{34}}{m_3}{m_4}{m_{34}} 
  \CG{j_{12}}{j_{34}}{j}{m_{12}}{m_{34}}{m}  
  =\\
  \sum_{j_{13},j_{24}}
  \alpha
  \NineJ{j_{1} }{j_{2} }{j_{12} }  
        {j_{3} }{j_{4} }{j_{34}}
        {j_{13}}{j_{24}}{j     }
  \CG{j_1}{j_3}{j_{13}}{m_1}{m_3}{m_{13}}
  \CG{j_{2}}{j_{4}}{j_{24}}{m_{2}}{m_{4}}{m_{24}} 
  \CG{j_{13}}{j_{24}}{j}{m_{13}}{m_{24}}{m}
  \,.
\end{multline}

\providecommand{\href}[2]{#2}\begingroup\raggedright\endgroup

\end{document}